\documentclass[reqno,11pt]{amsart}

\usepackage{amssymb}
\newtheorem{Proposition}{Proposition}[section]
\newtheorem{Definition}[Proposition]{Definition}
\newtheorem{Lemma}[Proposition]{Lemma}
\newtheorem{Theorem}[Proposition]{Theorem}
\newtheorem{Corollary}[Proposition]{Corollary}

\DeclareMathOperator{\Val}{Val}
\DeclareMathOperator{\Spin}{Spin}
\DeclareMathOperator{\G}{G_2}
\DeclareMathOperator{\SU}{SU}
\DeclareMathOperator{\SO}{SO}
\DeclareMathOperator{\U}{U}
\DeclareMathOperator{\Sp}{Sp}
\DeclareMathOperator{\vol}{vol}

\DeclareMathOperator{\Gr}{Gr}

\DeclareMathOperator{\Kl}{Kl}

\DeclareMathOperator{\Stab}{Stab}

\newcommand{\R}{\mathbb{R}}

\title{Integral geometry under $\G$ and $\Spin(7)$}
\author{Andreas Bernig}
 
\email{andreas.bernig@unifr.ch}
\address{D\'epartement de Math\'ematiques, Chemin du Mus\'ee 23, 1700 Fribourg, Switzerland}

\begin{document}

\begin{abstract}
A Hadwiger-type theorem for the exceptional Lie groups $\G$ and $\Spin(7)$ is proved. The algebras of $\G$ or $\Spin(7)$ invariant, translation invariant continuous valuations are both of dimension $10$. Geometrically meaningful bases are constructed and the algebra structures are computed. Finally, the kinematic formulas for these groups are determined.   
\end{abstract}

\thanks{ 
\\ Supported
  by the Schweizerischer Nationalfonds grant SNF PP002-114715/1.}
\maketitle 
\section{Introduction}

Let $V$ be a finite-dimensional vector space. By $\mathcal{K}(V)$ we denote the space of compact convex subsets of $V$. This space has a natural topology, the Hausdorff topology. A (convex) valuation is a functional $\mu:\mathcal{K}(V) \to \mathbb{C}$ with the following additivity property:
\begin{displaymath}
 \mu(K \cup L)=\mu(K)+\mu(L)-\mu(K \cap L), \quad \forall K,L, K \cup L \in \mathcal{K}(V).
\end{displaymath}
Examples of valuations are the Euler characteristic $\chi$ (which equals $1$ for non-empty $K$) and any Lebesgue measure.

In this paper, we will only consider continuous and translation invariant valuations. We set $\Val=\Val(V)$ for the vector space of all continuous and translation invariant valuations on $V$. Alesker \cite{ale04a}, \cite{ale05a} introduced a product structure on $\Val$ (in fact on some dense subspace consisting of {\it smooth valuations}). 

Let $V$ be endowed with a scalar product. For a compact subgroup $G$ of $\SO(V) \cong \SO(n)$, we let $\Val^G \subset \Val$ denote the subspace of $G$-invariant valuations. Alesker has shown that $\Val^G$ is finite-dimensional if and only if $G$ acts transitively on the unit sphere \cite{ale00}, \cite{ale06}. In this case, $\Val^G$ consists of smooth valuations and $\Val^G$ is a finite-dimensional algebra. 

Borel \cite{bor49}, \cite{bor50} and Montgomery-Samelson \cite{mosa43} gave a complete classification of connected compact Lie groups acting transitively and effectively on the sphere (compare also \cite{bes87}, 7.13). There are six infinite series 
\begin{equation} \label{eq_series}
 \SO(n), \U(n), \SU(n), \Sp(n), \Sp(n) \cdot \U(1), \Sp(n) \cdot \Sp(1) 
\end{equation}
and three exceptional groups 
\begin{equation} \label{eq_exceptions}
 \G, \Spin(7),\Spin(9). 
\end{equation}

For each group $G$ in this list, three fundamental problems of integral geometry are as follows. 
\begin{enumerate}
 \item[(i)] {\bf Hadwiger-type theorem for $G$}: determine the dimension of $\Val^G$ and give a geometrically meaningful basis of this space.
\item[(ii)] Compute the {\bf algebra structure} of $\Val^G$. 
\item[(iii)] Compute the whole array of {\bf kinematic formulas} and additive kinematic formulas.
\end{enumerate}

The classical case $G=\SO(n)$ was studied by Blaschke, Chern, Santal\'o, Hadwiger and others. Hadwiger showed that $\dim \Val^{\SO(n)}=n+1$, with a basis given by the intrinsic volumes $\mu_0=\chi,\ldots,\mu_n=\vol_n$. The graded algebra $\Val^{\SO(n)}$ is isomorphic to $\mathbb{C}[t]/(t^{n+1})$ (where $\deg t=1$). 

For a compact subgroup $G$ of $\SO(n)$, we set $\bar G:=G \ltimes V$ with the product of Haar and Lebesgue measure. The principal kinematic formula for $\SO(n)$ is
\begin{align*}
\int_{\overline{\SO(n)}} \chi(K \cap \bar g L)d\bar g =\sum_{k=0}^n \binom{n}{k}^{-1}\frac{\omega_k\omega_{n-k}}{\omega_n} \mu_k(K)\mu_{n-k}(L), \quad K,L \in \mathcal{K}(V). 
\end{align*}
Here and in the following, $\omega_n$ is the volume of the $n$-dimensional unit ball. Higher and additive kinematic formulas are also classical. 

Problem (i) for $G=\U(n)$ was solved by Alesker \cite{ale03b}, problem (ii) by Fu \cite{fu06}. In general, the product structure of $\Val^G$ determines all kinematic formulas \cite{befu06}, but it may be rather hard to write down explicit formulas. For $G=\U(n)$, this was recently achieved in \cite{befu08}, thus solving problem (iii). 

Alesker \cite{ale04} computed the dimension of $\Val^{\SU(2)}$. The algebra structure was computed in \cite{be06}. The kinematic formulas were known before by work of Tasaki \cite{tas03}. 

Problems (i)-(iii) for $G=\SU(n)$ were recently solved in \cite{be08} and we will use these results in an essential way in the present work.  

For the other groups from list \eqref{eq_series}-\eqref{eq_exceptions}, very few is known.

In this paper, we will solve problems (i)-(iii) for the exceptional Lie groups $\G$ and $\Spin(7)$. Our method uses the inclusions 
\begin{displaymath}
 \SU(3) < \G < \Spin(7), \SU(4) < \Spin(7)
\end{displaymath}
and the fact that the integral geometry of $\SU(3)$ and $\SU(4)$ is well-understood. In particular, we will need that $\G$ and $\Spin(7)$ act $2$-transitively on the unit sphere. This is not the case for $G=\Spin(9)$, which is the reason why other ideas have to be used in the study of the integral geometry under this group.


\section{Results}

We define the group $\G$ following \cite{br87}. Other references are \cite{hit00},  \cite{joy00} and \cite{fu_ha91}. A nice historical account can be found in \cite{agri07}.
 
Let $V$ be a $7$-dimensional real vector space. Given $\phi \in \Lambda^3 V^*$, we set 
\begin{displaymath}
 b(x,y):=\frac16 i_x \phi \wedge i_y \phi \wedge \phi \in \Lambda^7 V^*, \quad x,y \in V.
\end{displaymath}

If $b(x,x) \neq 0$ for all $ x \neq 0$, $\phi$ will be called {\it positive}. 

Choosing a basis $\omega_1,\ldots,\omega_7$ of $V^*$, we obtain an example of a positive $3$-form by setting 
\begin{align*}
\phi & = \omega_1 \wedge (\omega_2 \wedge \omega_3+\omega_4 \wedge \omega_5+\omega_6 \wedge \omega_7)\\
& \quad +\omega_2 \wedge \omega_4 \wedge \omega_6-\omega_2 \wedge \omega_5 \wedge \omega_7-\omega_3 \wedge \omega_4 \wedge \omega_7-\omega_3 \wedge \omega_5 \wedge \omega_6.
\end{align*}

\begin{Definition}
 Let $\phi \in \Lambda^3 V^*$ be positive. Then $\G=\G(V,\phi)$ is defined as 
\begin{displaymath}
 \G:=\{g \in GL(V)| g^*\phi=\phi\}.
\end{displaymath}
\end{Definition}

Since $\phi$ is positive, $b(x,x)$ lies in the same connected component of $\Lambda^7 V^* \setminus \{0\}$ for all $x \neq 0$. We fix the orientation on $V$ in such a way that $b(x,x)>0$ for all $x \neq 0$. 

Given an orientation preserving isomorphism $\tau:\Lambda^7 V^* \to \R$, we obtain a scalar product on $V$ by setting $\langle x,y\rangle:=\tau \circ b(x,y)$. This scalar product induces an isomorphism $\tau':\Lambda^7 V^* \to \R$, sending $\omega$ to 
 $\omega(e_1,\ldots,e_7)$ for a positively oriented orthonormal basis 
$e_1,\ldots,e_7$ of $V$. In general $\tau'\neq \tau$. However, there is a unique choice of $\tau$ such that $\tau'=\tau$. Indeed, rescaling $\tau$ by a positive factor $t^2$ will rescale $\tau'$ by $t^{-7}$, hence there is a unique $t$ with $\tau'=\tau$. We fix in the following the scalar product with this property. 

Since $b$ is $\G$-equivariant, $\G$ preserves $\langle \cdot,\cdot\rangle$ and the orientation. In other words, $\G \subset \SO(V)$. In particular, the intrinsic volumes $\mu_k,k=0,\ldots,7$ are $\G$-invariant.

Let $W \subset V$ be a $3$-dimensional subspace. Choose vectors $w_1,w_2,w_3 \in W$ such that $\|w_1 \wedge w_2 \wedge w_3\|=1$. We set $\phi(W):=\phi(w_1,w_2,w_3)$. Clearly, $\phi(W)$ is well-defined (i.e. independent of the choice of $w_1,w_2,w_3$) up to a sign; hence $\phi(W)^2$ is well-defined. It may be checked that two spaces $W_1,W_2 \in \Gr_3(V)$ are in the same $\G$-orbit if and only if $\phi(W_1)^2=\phi(W_2)^2$, but we will not use this fact in the following. 

For a face $F$ of a convex polytope $P \subset V$, we let $W_F$ be the linear space parallel to $F$. The normalized volume of the outer angle at $F$ is denoted by $\gamma(F)$. 

\begin{Theorem} (Hadwiger-type theorem for $\G$)\\ \label{thm_hadwiger_g2}
 There is a unique valuation $\nu_3 \in \Val^{\G,+}$ such that for each polytope $P \subset V$ 
\begin{equation} \label{eq_def_tau3}
 \nu_3(P) =\sum_{F,\dim F=3} \gamma(F) \phi(W_F)^2 \vol(F).
\end{equation}
There is a unique valuation $\nu_4 \in \Val^{\G,+}$ such that for each polytope $P \subset V$ 
\begin{equation} \label{eq_def_tau4}
 \nu_4(P) =\sum_{F,\dim F=4} \gamma(F) \phi(W_F^\perp)^2 \vol(F).
\end{equation}
The valuations $\mu_0,\ldots,\mu_7,\nu_3,\nu_4$ form a basis of the vector space $\Val^{\G}$. In particular, 
\begin{displaymath}
 \dim \Val^{\G}=10;
\end{displaymath}
and all $\G$-invariant valuations are even.   
\end{Theorem}

\begin{Corollary} \label{cor_only_even}
 Let $G$ be any of the groups from the list \eqref{eq_series}-\eqref{eq_exceptions}. Then $\Val^G$ consists only of even valuations. 
\end{Corollary}

\proof
If $-1 \in G$, then the result is clear. The only groups with $-1 \notin G$ are $SO(n)$ and $SU(n)$ for $n$ odd and $G_2$. The case $G=SO(n)$ is the classical Hadwiger theorem, which shows in particular that $SO(n)$-invariant valuations are even. The geometric reason behind this is a lemma of Sah (\cite{klro97}, Prop. 8.3.1).

The case $G=SU(n)$ was studied in \cite{be08}, again there are no odd invariant valuations. Finally, the case $G=G_2$ is treated in Theorem \ref{thm_hadwiger_g2}.
\endproof

We set 
\begin{align*}
 \nu_3' & :=5\nu_3-\mu_3\\
\nu_4' & := 5 \nu_4-\mu_4.
\end{align*}

\begin{Theorem} (Algebra structure of $\Val^{\G}$)\\ \label{thm_algebra_g2}
 Let $t$ and $u$ be variables of degree $1$ respectively $3$. Then the map $t \mapsto \frac{2}{\pi}\mu_1, u \mapsto \frac{\nu_3'}{\pi^2}$ covers an isomorphism between graded algebras 
\begin{displaymath}
 \Val^{\G} \cong \mathbb{C}[t,u]/(t^8,t^2u,u^2+4t^6). 
\end{displaymath}
\end{Theorem}

By the methods of \cite{befu06}, one can translate this algebra structure into the whole array of kinematic formulas for $\G$. As an example, we write down in explicit form the principal kinematic formula for $\G$. 

\begin{Corollary} (Principal kinematic formula for $\G$)\\  \label{cor_kin_g2}
For compact convex sets $K,L \subset V$ we have  
\begin{align*}
\int_{\bar \G} \chi(K \cap \bar g L)d\bar g 
=\int_{\overline{\SO(V)}} \chi(K \cap \bar g L)d\bar g + \frac{1}{2^9} \nu_3'(K)\nu_4'(L)+\frac{1}{2^9} \nu_4'(K)\nu_3'(L).
\end{align*}
\end{Corollary}

Now we turn our attention to $\Spin(7)$. We can define it in a short way as the universal covering of $\SO(7)$. However, it will be convenient to have an explicit description similar to the one for $\G$ given above. 

Let $V_+$ be a complex $4$-dimensional hermitian vector space. Let $\Omega$ be the symplectic $2$-form of $V_+$ and $\beta \in \Lambda^4 V_+^*$ be the real part of the holomorphic volume form, i.e. $\beta(v_1,v_2,v_3,v_4)=Re \det(v_1,v_2,v_3,v_4)$.  Set 
\begin{equation} \label{eq_def_Phi}
 \Phi:=\frac12 \Omega \wedge \Omega + \beta \in \Lambda^4V^*. 
\end{equation}

\begin{Definition}
\begin{displaymath}
 \Spin(V_+):=\left\{ g \in Gl(V_+)| g^*\Phi=\Phi\right\}.
\end{displaymath}
\end{Definition}

This group is a compact connected subgroup of $\SO(V_+)$ which acts transitively on the unit sphere of $V_+$. For $v \in S(V_+)$, the stabilizer of $\Spin(V_+)$ at $v$ is $\G(W,\phi)$, where $W:=T_vS(V_+)$ and 
\begin{equation} \label{eq_def_phi}
 \phi:=*_W( \Phi|_W) \in \Lambda^3W^*. 
\end{equation}

For $W \in \Gr_4(V_+)$ we set $\Phi(W):=\Phi(w_1,w_2,w_3,w_4)$, where $w_1,w_2,w_3,w_4 \in W$ satisfy $\|w_1 \wedge w_2 \wedge w_3 \wedge w_4\|=1$. Then $\Phi(W)$ is well-defined up to a sign. It is easy to show that $\Spin(V_+)$ acts transitively on $\Gr_k(V_+), k \neq 4$ and that $W_1,W_2 \in \Gr_4(V_+)$ belong to the same $\Spin(V_+)$-orbit if and only if $\Phi(W_1)^2=\Phi(W_2)^2$.   
 
\begin{Theorem} (Hadwiger-type theorem for $\Spin(7)$) \\
\label{thm_hadwiger_spin7}
 There exists a unique valuation $\eta \in \Val^{\Spin(V_+)}$ such that for each polytope $P \subset V_+$ we have 
\begin{equation} \label{eq_def_eta}
 \eta(P)=\sum_{F,\dim F=4} \gamma(F) \vol(F) \Phi(W_F)^2. 
\end{equation}
The valuations $\mu_0,\ldots,\mu_8,\eta$ form a basis of $\Val^{\Spin(V_+)}$, in particular 
\begin{displaymath}
\dim \Val^{\Spin(V_+)}=10. 
\end{displaymath}
\end{Theorem}

We set $\eta':=5\eta-\mu_4$. It will turn out that $\eta'$ is {\it primitive} in the sense that $\mu_1 \cdot \eta'=0$.  

\begin{Theorem} (Algebra structure of $\Val^{\Spin(7)}$) \\
\label{thm_algebra_spin7}
Let $V_+$ be, as before, a hermitian vector space of complex dimension $4$. Let $t$ and $u$ be variables of degree $1$ respectively $4$. Then the map $t \mapsto \frac{2}{\pi}\mu_1, u \mapsto \frac{\eta'}{\pi^2}$ covers an isomorphism between graded algebras 
\begin{displaymath}
 \Val^{\Spin(V_+)} \cong \mathbb{C}[t,u]/(t^9,u^2-t^8,ut). 
\end{displaymath}
\end{Theorem}

Again, we can compute all kinematic formulas for $\Spin(V_+)$ using this theorem and the results of \cite{befu06}. We only give the following example. 

\begin{Corollary} (Principal kinematic formula for $\Spin(7)$)\\
For compact convex sets $K,L \subset V_+$ we have  
\begin{displaymath}
\int_{\overline{\Spin(V_+)}} \chi(K \cap \bar g L)d\bar g=\int_{\overline{\SO(V_+)}} \chi(K \cap \bar g L)d\bar g+\frac{3}{7!}\eta'(K)\eta'(L).
\end{displaymath}
\end{Corollary}


\section{Tools} 
\label{sec_tools}

Let $V$ be a vector space of dimension $n$. Let $\Val:=\Val(V)$ be the vector space of continuous, translation invariant valuations. A valuation $\mu \in \Val$ has degree $k$ if $\mu(tK)=t^k \mu(K)$ for all $t \geq 0$. It is even (resp. odd) if $\mu(-K)=\mu(K)$ (resp. $\mu(-K)=-\mu(K)$). We write $\Val_k^\pm$ for the subspace of valuations of degree $k$ and parity $\pm$. 
\begin{Theorem} (P. McMullen, \cite{mcmu77})\\
There is a direct sum decomposition
 \begin{displaymath}
 \Val=\bigoplus_{\substack{k=0,\ldots,n\\ \epsilon=\pm}} \Val_k^\epsilon.
\end{displaymath}

\end{Theorem}
 
Suppose that $V$ is endowed with a Euclidean scalar product. Let $G$ be a compact subgroup of $\SO(V)$ acting transitively on the unit sphere. Alesker has shown that $\dim \Val^G<\infty$. He introduced on $\Val^G$ a structure of a commutative and associative algebra with unit $\chi$ satisfying Poincar\'e duality \cite{ale04a,ale05a}. 
 
The next theorem is a particular case of the more general hard Lefschetz theorem of \cite{bebr07} (where $G$-invariance is replaced by smoothness in the sense of \cite{ale05a}). 

\begin{Theorem} (Hard Lefschetz theorem, \cite{bebr07})\\ \label{thm_hlt}
 Let $\tilde \Lambda:\Val^G_*(V) \to \Val^G_{*-1}(V)$ be defined by
\begin{displaymath}
 \tilde \Lambda \mu(K):=\left.\frac{d}{dt}\right|_{t=0} \mu(K+tB).
\end{displaymath}
Then $\tilde \Lambda^{2k-n}:\Val_k^G(V) \to \Val_{n-k}^G(V)$ is an isomorphism. 
In particular, the {\it Betti numbers} 
\begin{displaymath}
 h_k:=\dim \Val_k^G
\end{displaymath}
satisfy the inequalities
\begin{align*}
 h_0 & \leq h_1 \leq \ldots \leq h_{\lfloor n/2\rfloor}\\
 h_{\lfloor n/2\rfloor} & \geq h_{\lfloor n/2\rfloor+1} \geq \ldots \geq h_n
\end{align*}
and the equations 
\begin{displaymath}
 h_k=h_{n-k}, \quad k=0,\ldots,n. 
\end{displaymath}
\end{Theorem}

The theorem in the case of {\it even} valuations was proved before by Alesker \cite{ale03b}, using the following fundamental result.  

\begin{Theorem} (Klain, \cite{kl95})\\ \label{thm_klain}
 Let $\mu \in \Val^+(V)$ be simple, i.e. $\mu$ vanishes on lower-dimensional compact convex bodies. Then $\mu$ is a multiple of the Lebesgue measure on $V$, in particular $\mu$ is of degree $n$.  
\end{Theorem}

This theorem has the following important consequence. The restriction of $\mu \in \Val_k^+$ to a subspace $W \in \Gr_k(V)$ is a multiple of the Lebesgue measure on $W$. Putting $\Kl_\mu$ for the proportionality factor, one gets a function $\Kl_\mu \in C(\Gr_k V)$, called the {\bf Klain function} of $\mu$. The resulting map $\Val_k^+ \to C(\Gr_k V)$ is injective. 

If $\mu$ is smooth in the sense of \cite{ale05a} (in particular if $\mu \in \Val^G$ for one of the groups from \eqref{eq_series}-\eqref{eq_exceptions}), then there exists a unique valuation $\hat \mu \in \Val_{n-k}^+$ whose Klain function is given by $\Kl_{\hat \mu}=\Kl_\mu \circ \perp$. This valuation is called {\bf Akesker dual} or {\bf Fourier transform} of $\mu$. The Fourier transform can be extended to odd valuations, but we will not need this here.

An analogue of Klain's injectivity result in the case of {\it odd} valuations was proved by Schneider. 
\begin{Theorem} (Schneider, \cite{schn96})\\  \label{thm_schn}
Let $\mu \in \Val^-(V)$ be simple. Then there exists an odd continuous function $f$ on the unit sphere $S(V)$ such that 
\begin{displaymath}
 \mu(K)=\int_{S(V)} f(v) dS_{n-1}(K,v),
\end{displaymath}
where $S_{n-1}(K,\cdot)$ is the $n-1$-th surface area measure of $K$. In particular, $\mu$ is of degree $n-1$. 
\end{Theorem}

\begin{Corollary} \label{cor_dimension_bounds}
Let $V$ be a Euclidean vector space of dimension $n$. Let $G<\SO(V)$ act transitively on the unit sphere. For a hyperplane $W \subset V$, let $H:=\Stab_G(W)$. Consider the restriction map 
\begin{align*}
 r_k^{\pm}: \Val^{G,\pm}_k(V) & \to \Val^{H,\pm}_k(W)\\
\mu & \mapsto \mu|_W.
\end{align*}
Then $r_k^+$ is injective for $k \neq n$ while $r_k^-$ is injective for $k \neq n-1$.   
\end{Corollary}

\proof
Since $G$ acts transitively on $\Gr_{n-1}(V)$, a valuation $\mu \in \ker r_k^\pm$ is simple. The statement thus follows from Theorems \ref{thm_klain} and \ref{thm_schn}.  
\endproof

If $H$ acts transitively on the unit sphere of $W$, then the dimension of $\Val_k^{H,\pm}$ is finite and we get a bound for the dimension of $\Val_k^{G,\pm}$. This is the case precisely when $G$ acts $2$-transitively on the unit sphere $S(V)$ (this means that if $v_1,v_2,v_1',v_2' \in S(V)$ with $\langle v_1,v_2\rangle=\langle v_1',v_2'\rangle$, then there exists $g \in G$ with $gv_1=v_1'$ and $gv_2=v_2'$). For the groups $G$ from the list \eqref{eq_series}-\eqref{eq_exceptions}, only $\SO(n)$, $\G$ and $\Spin(7)$ have this property. 

We will need some results from \cite{befu08} and \cite{be08}. Let $V$ be a hermitian vector space of complex dimension $n$ and $W \in \Gr_k(K)$. Set $p:=\lfloor k/2\rfloor$. The restriction of the symplectic form $\Omega$ of $V$ to $W$ can be written in the form
\begin{displaymath}
 \Omega|_W=\sum_{i=1}^p \cos(\theta_i) \alpha_{2i-1} \wedge \alpha_{2i}.
\end{displaymath}
  Here $\alpha_1,\ldots,\alpha_k$ is an orthonormal basis of the dual space $W^*$ and $0 \leq \theta_1 \leq \ldots \leq \theta_p \leq \frac{\pi}{2}$ is called the {\bf multiple K\"ahler angle} of $W$ \cite{tas03b}. We denote the $q$-th elementary symmetric function by $\sigma_q$. 

\begin{Theorem} (\cite{befu08})\\ \label{thm_tasaki_valuations}
For $0 \leq q \leq p$, there exists a unique valuation $\tau_{k,q} \in \Val_k^{\U(V)}$ such that for a polytope $P \subset V$ 
\begin{displaymath}
 \tau_{k,q}(P)=\sum_{F, \dim F=k} \gamma(F)\vol(F) \sigma_q(\cos^2 \theta_1(W_F),\ldots,\cos^2 \theta_p(W_F)).
\end{displaymath}
\end{Theorem}
These valuations are called {\bf Tasaki valuations}. Note that $\tau_{k,0}=\mu_k$, the $k$-th intrinsic volume. 

Under the subgroup $\SU(V)<\U(V)$, there are some more invariant valuations; they are even and appear in the middle degree $n$. More precisely, 
\begin{align*}
\Val_k^{\SU(V)} & =\Val_k^{\U(V)} \quad \text{ if } k \neq n;\\
\dim \Val_n^{\SU(V)} & = \Val_k^{\U(V)}+4 \quad \text{ if } n \equiv 0 \mod 2;\\
\dim \Val_n^{\SU(V)} & = \Val_k^{\U(V)}+2 \quad \text{ if } n \equiv 1 \mod 2.
\end{align*}

Let us describe the ``new'' invariant valuations more explicitly. 

For $W \in \Gr_n(V)$, we set $\Theta(W):=\det(w_1,\ldots,w_n)$, where $w_1,\ldots,w_n$ is a basis of $W$ with $\|w_1 \wedge \ldots \wedge w_n\|=1$. If $\Omega|_W$ is non-degenerated (which happens if and only if $n=2p$ is even and $ \theta_p(W) < \frac{\pi}{2}$), we ask in addition that $w_1,\ldots,w_n$ be a positively oriented basis. 

If $\Omega|_W$ is non-degenerated, then $\Theta(W) \in \mathbb{C}$ is a well-defined invariant. Otherwise, $\Theta(W) \in \mathbb{C}/\{\pm 1\}$ is well-defined. It is easily checked (compare \cite{be08}) that in both cases
\begin{equation} \label{eq_norm_Theta}
|\Theta(W)|=\prod_{i=1}^p \sin \theta_i(W).  
\end{equation}

\begin{Theorem} (\cite{be08})\\ \label{thm_sun_valuations}
There exists a unique valuation $\phi_{n,2} \in \Val_n^{\SU(V)}$ such that 
\begin{displaymath}
 \phi_{n,2}(P) =\sum_{F, \dim F=n} \vol(F) \gamma(F) \Theta(W_F)^2.
\end{displaymath}
If $n=2p$ is even, there exists a unique valuation $\phi_{n,1} \in \Val_n^{\SU(V)}$ such that 
\begin{displaymath}
\phi_{n,1}(P) =\sum_{F, \dim F=n} \vol(F) \gamma(F) \Theta(W_F) \prod_{j=1}^p \cos \theta_j(W_F).
\end{displaymath}
\end{Theorem}
Note that these valuations are complex-valued. They were denoted by $\phi_1,\phi_2$ in \cite{be08}. 


\section{Construction of invariant valuations}

In this section, we let $V_+$ be a complex $4$-dimensional hermitian vector space with the $4$-form $\Phi$ as in \eqref{eq_def_Phi} and $V \subset V_+$ a hyperplane with the positive $3$-form $\phi$ from \eqref{eq_def_phi}. We will write $\G$ instead of $\G(V,\phi)$.  

By Corollary \ref{cor_dimension_bounds} with $G:=\G$ and $H:=\SU(3)$, we get the following dimensions: 

\begin{align*}
 \dim \Val^{\G,-} & =0\\
 \dim \Val_k^{\G,+} & =1 \quad k \neq 3,4\\
 \dim \Val_3^{\G,+} & = \dim \Val_4^{\G,+}\leq 2. 
\end{align*}

In particular, there are no non-zero odd $\G$-invariant valuations, which is non-trivial since $-1 \notin \G$.  

Clearly, there are no odd $\Spin(V_+)$-invariant valuations, since $-1 \in \Spin(V_+)$. We apply Corollary \ref{cor_dimension_bounds} with $G:=\Spin(V_+)$ and $H:=\G$ and find the following dimensions: 
\begin{align*}
 \dim \Val_k^{\Spin(V_+)} & = 1, \quad k \neq 4\\
\dim \Val_4^{\Spin(V_+)} &  \leq 2. 
\end{align*}

Let us show that in fact $\dim \Val_4^{\Spin(V_+)}=2$, which also implies that $\dim \Val_3^{\G}=\dim \Val_4^{\G}=2$. 

Recall that $\SU(V_+)$ is a subgroup of $\Spin(V_+)$. Besides the intrinsic volume $\mu_4$, there is another $\SU(V_+)$-invariant valuation of degree $4$ which remains invariant under the bigger group $\Spin(V_+)$.

\begin{Proposition}
$\eta:=\frac12 \tau_{4,0}-\frac12 \tau_{4,1}+\frac32 \tau_{4,2}+\frac12 Re \phi_{4,2}+2 Re \phi_{4,1} \in \Val_4^{\SU(V_+)}$ is $\Spin(V_+)$-invariant and satisfies \eqref{eq_def_eta}. 
\end{Proposition}

\proof
Let $W \in \Gr_4(V)$. Fix a basis $w_1,\ldots,w_4$ of $W$ with $\|w_1 \wedge \ldots \wedge w_4\|=1$. If $\Omega|_W$ is non-degenerated, we ask in addition that $w_1 \wedge \ldots \wedge w_4$ be positive. 

Let us write $\theta_i:=\theta_i(W), i=1,2$; $\sigma_q:=\sigma_q(\cos \theta_1,\cos \theta_2), q=0,1,2$ and $\Theta:=\Theta(W)$. Then, by \eqref{eq_def_Phi} and \eqref{eq_norm_Theta}, 
\begin{align*}
 \Phi(W)^2 & = \left(\cos \theta_1 \cos \theta_2+Re \Theta\right)^2 \\
& = \cos^2 \theta_1 \cos^2 \theta_2+(Re \Theta)^2+2 \cos \theta_1 \cos \theta_2 Re \Theta\\
& = Re\left(\cos^2 \theta_1 \cos^2 \theta_2+\frac12 \Theta^2+\frac12 sin^2 \theta_1 \sin^2 \theta_2+2 \cos \theta_1 \cos \theta_2 \Theta\right)\\
& = \frac12 \sigma_0-\frac12 \sigma_1 +\frac32 \sigma_2+Re\left(\frac12 \Theta^2+2 \cos \theta_1 \cos \theta_2 \Theta\right).
\end{align*}

Using Theorems \ref{thm_tasaki_valuations} and \ref{thm_sun_valuations} we get for each polytope $P \subset V_+$  
\begin{displaymath}
 \eta(P)=\sum_{F, \dim F=4} \vol(F) \gamma(F) \Phi(W_F),
\end{displaymath}
which is \eqref{eq_def_eta}. This equation shows in particular that $\eta$ is $\Spin(V_+)$ invariant. 
\endproof

Theorem \ref{thm_hadwiger_spin7} follows from the proposition. 

Let $V$ be a hyperplane in $V_+$. The stabilizer of $\Spin(V_+)$ at $V$ is $\G=\G(V,\phi)$, where $\phi:=*_W \Phi|_W$ (note that $\phi$ is only well-defined up to a sign depending on the choice of orientation on $V$). 

We set $\nu_4:=\eta|_V \in \Val_4^{\G}$. From \eqref{eq_def_phi} we deduce that  
\begin{displaymath}
 \phi(W^\perp)^2=\Phi(W)^2, \quad W \in \Gr_4(V).
\end{displaymath}
Equation \eqref{eq_def_tau4} follows immediately.

We define $\nu_3$ by 
\begin{displaymath}
 \nu_3 := \hat \nu_4 \in \Val_3^{\G}.
\end{displaymath}

Unfortunately, the verification of \eqref{eq_def_tau3} seems to be less simple. If a continuous translation invariant valuation satisfies \eqref{eq_def_tau3}, then clearly its Klain function is given by $W \mapsto \phi(W)^2$, which is also the Klain function of $\nu_3=\hat \nu_4$. Without giving the details, we indicate that \eqref{eq_def_tau3} follows by noting that $\nu_3$ is a constant coefficient valuation in the sense of \cite{befu08}. In fact, if $N_1(K) \in \mathcal{I}_7(V \times V)$ denotes the disk bundle of $K$, then we can represent $\nu_3$ by 
\begin{displaymath}
 \nu_3(K)=\frac{1}{2\pi^2} N_1(K)(p_1^*\phi \wedge p_2^* *\phi).
\end{displaymath}
Here $p_1,p_2:V \times V \to V$ are the natural projections. Note that $p_1^*\phi \wedge p_2^* *\phi \in \Lambda^7 (V^* \times V^*)$ has constant coefficients.


\section{Algebra structure and kinematic formulas}
 
\proof[Proof of Theorem \ref{thm_algebra_spin7}]
Since $\SU(V_+) \subset \Spin(V_+)$, we get a canonical injection of graded algebras
\begin{displaymath}
 \Val^{\Spin(V_+)} \hookrightarrow \Val^{\SU(V_+)}.
\end{displaymath}
The image of $\eta$ in $\Val^{\SU(V_+)}$ is $\frac12 \tau_{4,0}-\frac12 \tau_{4,1}+\frac32 \tau_{4,2}+\frac12 Re \phi_{4,2}+2Re \phi_{4,1}$. The valuations $\phi_{4,2}$ and $\phi_{4,1}$ are in the annihilator of $\mu_1$. Using the results of \cite{befu08} and \cite{be08}, one computes that 
\begin{align*}
 \mu_1 \cdot \eta  & = \mu_1 \cdot \left(\frac12 \tau_{4,0}-\frac12 \tau_{4,1}+\frac32 \tau_{4,2}\right)= \frac{8}{15} \mu_5,\\
\eta \cdot \eta & = \left(\frac12 \tau_{4,0}-\frac12 \tau_{4,1}+\frac32 \tau_{4,2}+\frac14 \phi_{4,2}+\frac14 \bar \phi_{4,2}+\phi_{4,1}+\bar \phi_{4,1}\right)^2\\
& = \left(\frac12 \tau_{4,0}-\frac12 \tau_{4,1}+\frac32 \tau_{4,2}\right)^2+\frac{1}{8} \phi_{4,2}\bar \phi_{4,2}+2\phi_{4,1}\bar \phi_{4,1}\\
& = \frac{203}{3} \mu_8.
\end{align*} 

Setting $\eta':=5\eta-\mu_4$ one gets $\eta' \cdot \mu_1=0$ and $\eta' \cdot \eta'=\frac{7!}{3} \mu_8$.  
\endproof

Next, we study the algebra structure of $\Val^{\G}$. Fix a hyperplane $V_0 \subset V$. Let $x \in V_0^\perp$ be of norm $1$. Then $i_x \phi$ is a symplectic form on $V_0$ and there exists a unique almost complex structure $J$ on $V_0$ with $\langle v,w\rangle=i_x \phi(Jv,w)$. 

The stabilizer of $\G$ at $V_0$ is given by $\SU(V_0) \cong \SU(3)$. 

\begin{Lemma}
Let $r:\Val^{\G} \to \Val^{\SU(V_0)}$ be the natural restriction map. Then 
\begin{align*}
 r(\nu_3) & =\frac12 \tau_{3,0}-\frac12 \tau_{3,1}+\frac12 Re \phi_{3,2}\\
r(\nu_4) & = \tau_{4,2}.
\end{align*}
\end{Lemma}

\proof
For $W \in \Gr_3(V_0)$, set $\theta:=\theta_1(W)$ and $\Theta:=\Theta(W)$. Then $\phi(W)=Re \Theta(W)$ and therefore 
\begin{displaymath}
 \phi(W)^2=(Re \Theta)^2=\frac12 Re(\Theta^2+|\Theta|^2)=\frac12 Re(\Theta^2+1-\cos^2 \theta).
\end{displaymath}

For $W \in \Gr_4(V_0)$, we have $\phi(W^\perp)=\cos \theta_2=\sigma_2(\cos \theta_1, \cos \theta_2)$, where $0=\theta_1 \leq \theta_2$ are the K\"ahler angles of $W$. 

The statement of the lemma now follows from the definitions in Theorem \ref{thm_tasaki_valuations} and Theorem \ref{thm_sun_valuations}.
\endproof

\proof[Proof of Theorem \ref{thm_algebra_g2}]
The Alesker product is compatible with the restriction $r:\Val^{\G} \to \Val^{\SU(V_0)}, \mu \mapsto \mu|_{V_0}$, hence 
\begin{displaymath}
 r(\mu_2 \cdot \nu_3) =\frac12 \mu_2 \cdot (\tau_{3,0}-\tau_{3,1})=\frac{4}{5} \mu_5. 
\end{displaymath}

Similarly, 
\begin{align*}
 r(\nu_3 \cdot \nu_3) & =\left(\frac12 \tau_{3,0}-\frac12 \tau_{3,1}+\frac14 \phi_{3,2}+\frac14 \bar \phi_{3,2}\right)^2\\
& =\frac14 (\tau_{3,0}-\tau_{3,1})^2+\frac{1}{8} \phi_{3,2} \bar \phi_{3,2}\\
& = -\frac{9\pi}{8} \mu_6.
\end{align*}
Since $r_5^+$ and $r_6^+$ are injective, we deduce that $\mu_{2} \cdot \nu_3=\frac45 \mu_5$ and $\nu_3 \cdot \nu_3=-\frac{9\pi}{8} \mu_6$ in $\Val^{\G}$. From these two equations, the theorem follows. 
\endproof

\proof[Proof of Corollary \ref{cor_kin_g2}]
Let us write $\nu_4$ in terms of the basis $\mu_4,\mu_1 \cdot \nu_3$ of $\Val_4^{\G}$. As above, we get that $r(\mu_1 \cdot \nu_3)=\frac{3}{16}\pi (\mu_4-\tau_{4,2})$ and therefore
\begin{displaymath}
 \nu_4=\mu_4-\frac{16}{3\pi}\mu_1 \cdot \nu_3.
\end{displaymath}
Corollary \ref{cor_kin_g2} now follows from the methods in \cite{befu06}. 
\endproof
  
We rescale the derivation operator $\tilde \Lambda$ of Theorem \ref{thm_hlt} as in \cite{befu08} by setting, on $\Val_k^{G_2}$,
\begin{displaymath}
 \Lambda:=\frac{\omega_{7-k}}{\omega_{8-k}} \tilde \Lambda
\end{displaymath}
and rescale multiplication by $\mu_1$ by 
\begin{displaymath}
 L \mu:=\frac{2\omega_k}{\omega_{k+1}} \mu_1 \cdot \mu.
\end{displaymath}

\begin{Corollary}
 The operators $L,\Lambda$ define, together with the degree counting operator $H\mu=(2\deg \mu-7) \mu$ a representation of the Lie algebra $\mathfrak{sl}_2$ on $\Val^{\G}$.
\end{Corollary}

\proof
The commutator relations $[H,L]=2L$ and $[H,\Lambda]=-2\Lambda$ are immediate. We have to show that $[L,\Lambda] \mu =H\mu$ for each $\mu \in \Val^{\G}$. If $\mu$ is one of the intrinsic volumes, then this follows from a direct computation. It remains to show that this equation also holds for $\mu=\nu_3$ and $\mu=\nu_4$. The Fourier transform intertwines $L$ and $\Lambda$ (compare \cite{befu06}) and therefore
\begin{displaymath}
 \Lambda \nu_3=\Lambda \hat \nu_4=\widehat{L \nu_4}=\hat \mu_5=\mu_2. 
\end{displaymath}
Similarly,
\begin{displaymath}
 \Lambda \nu_4=\Lambda \hat \nu_3=\widehat{L \nu_3}=\hat \mu_4-\hat \nu_4=\mu_3-\nu_3. 
\end{displaymath}
It is now easy to compute that $[L,\Lambda]\nu_3=-\nu_3$ and $[L,\Lambda]\nu_4=\nu_4$.
\endproof

Let us conclude with some general remarks. 
The algebra structure of the spaces $\Val^G$ is now known for the groups $G=\SO(n),\U(n),\SU(n),\G,\Spin(7)$ from the list \eqref{eq_series}-\eqref{eq_exceptions} . By looking at these algebras, some natural questions arise. 
\begin{itemize}
 \item For $G$ as above, rescaled versions of the derivation operator and of multiplication by $\mu_1$ define an $\mathfrak{sl}_2$-representation on $\Val^G$. We conjecture that this will also be the case for the remaining groups from the list  \eqref{eq_series}-\eqref{eq_exceptions}. On the bigger space $\Val^{sm}$ of smooth translation invariant valuations, this is not true. 
\item By Corollary \ref{cor_only_even}, $\Val^G \subset \Val^+$ for each of these groups $G$. However, our proof relies on a non-geometric case by case argument. Is there a direct explanation of this fact?
\item For each such group $G$, $\Val^G$ consists of constant coefficient valuations in the sense of \cite{befu08}. Is this also true for the remaining groups from the list \eqref{eq_series}-\eqref{eq_exceptions}? 
\end{itemize}


\end{document}